\documentclass[11pt]{amsart}
\usepackage{amsmath, amssymb, amsthm}
\usepackage{hyperref}
\usepackage{tikz}
\usepackage{subfigure}
\usepackage{ytableau}
\usepackage{array}
\usepackage{amsfonts}
\usepackage{mathtools}
\theoremstyle{plain}
\newtheorem{theorem}{Theorem}[section]
\newtheorem{lemma}[theorem]{Lemma}
\newtheorem{corollary}[theorem]{Corollary}
\newtheorem{proposition}[theorem]{Proposition}

\newtheorem{definition}{Definition}
\newtheorem{example}[]{Example}
\newtheorem{remark}[]{Remark}

\newtheorem{case}[]{Case}

\renewcommand{\t}{\tau}

\newcommand{\cP}{\mathcal{P}}
\newcommand{\bs}{\mathrm{bs}}
\numberwithin{equation}{section}

\begin{document}
	\title[Divisibility of Character Values of Wreath Products of Symmetric Groups]
	{Divisibility of Character Values of Wreath Products of Symmetric Groups}
	\author{JAYANTA BARMAN}
	\address{JAYANTA BARMAN\\ Indian Institute of Technology Kharagpur \\Department of Mathematics \\
		Kharagpur-721302, India.} 
	\email{b1999jayanta@gmail.com, b1999jayanta@kgpian.iitkgp.ac.in}
	
	\author[Kamalakshya Mahatab]{Kamalakshya Mahatab}
	\address{Kamalakshya Mahatab\\ Indian Institute of Technology Kharagpur \\Department of Mathematics \\
		Kharagpur-721302, India.} 
	\email{accessing.infinity@gmail.com, kamalakshya@maths.iitkgp.ac.in}
	
	
	\begin{abstract}
		Let $G$ be a finite group whose character table is
		integer-valued, and let $K$ be any fixed positive
		integer. We prove that almost every entry of the
		character table of the wreath product $G\wr S_N$
		is divisible by $K$ as $N\to\infty$.
	\end{abstract}

	\maketitle
	
	\section{Introduction}\label{section1}
	
	Let $S_N$ be the symmetric group on $N$ letters, and let $G$ be a
	finite group with $k$ conjugacy classes. By the Murnaghan--Nakayama
	rule~\cite[Theorem~2.4.7]{james}, the irreducible characters of
	$S_N$ are integer-valued. Moreover, the number of irreducible
	characters of $S_N$ is equal to $p(N)$, the number of partitions
	of $N$.
	
	Based on extensive numerical calculations, Miller
	\cite{miller1,miller2} conjectured that, for every fixed prime
	power $q$, the proportion of entries in the character table of
	$S_N$ divisible by $q$ tends to $1$ as $N\to\infty$.
	Peluse and Soundararajan \cite{peluse3} proved this conjecture.
	They also proved that, for every fixed positive integer $K$,
	almost all entries in the character table of $S_N$ are divisible
	by $K$ as $N\to\infty$.
	
	These results naturally lead to the question of whether similar
	divisibility properties hold for other families of groups with
	integer-valued character tables. In this paper, we study this
	question for the wreath product $G\wr S_N$. By the
	Murnaghan--Nakayama rule~\cite[Theorem~4.4.10]{james}
	(see also \cite[Proposition~2.8]{dong}), if $G$ has an
	integer-valued character table, then the character table of
	$G\wr S_N$ is also integer-valued. Thus, it is natural to
	investigate the divisibility properties of its character values. We also mention that, in our recent work \cite{barman4},
	we obtained a lower bound for the number of zeros in
	the character table of $G\wr S_N$ without assuming
	that $G$ has an integer-valued character table.
	
	Dong et al.~\cite{dong} proved that, for every fixed prime $p$,
	almost all entries in the character table of $G\wr S_N$ are
	divisible by $p$. Our first result extends their theorem from
	primes to prime powers. It also generalizes the result of
	Peluse and Soundararajan \cite{peluse3} from symmetric groups
	to wreath products.
	
	\begin{theorem}\label{thm-main}
		Let $G$ be a finite group with an integer-valued character
		table, and let $G\wr S_N$ be the wreath product of $G$ with
		the symmetric group $S_N$. Then, for every fixed prime power
		$q$, the proportion of entries in the character table of
		$G\wr S_N$ divisible by $q$ tends to $1$ as $N\to\infty$.
	\end{theorem}
	
	The proof of Theorem~\ref{thm-main} uses the combinatorial
	description of the irreducible representations and conjugacy
	classes of $G\wr S_N$. If $G$ has $k$ conjugacy classes,
	then both the irreducible representations and the conjugacy
	classes of $G\wr S_N$ are indexed by $k$-multipartitions of $N$.
	
	A $k$-multipartition of a positive integer $N$ is a sequence
	$\lambda=(\lambda^{(1)},\lambda^{(2)},\ldots,\lambda^{(k)})$,
	where each $\lambda^{(j)}$ is an integer partition, possibly
	empty, satisfying
	$\sum_{j=1}^{k}|\lambda^{(j)}|=N$.
	Here, $|\lambda^{(j)}|$ denotes the sum of the parts of
	$\lambda^{(j)}$.
	
	James and Kerber \cite[Corollary~4.4.4]{james} showed that
	the number of irreducible representations of $G\wr S_N$
	is given by
	\begin{equation*}
		\sum_{\substack{N=n_1+n_2+\cdots+n_k\\n_j\geq 0}}
		p(n_1)p(n_2)\cdots p(n_k).
	\end{equation*}
	This sum is precisely the number of $k$-multipartitions of $N$. Counting $k$-multipartitions according to the
	size of their first component, we obtain
	\begin{equation}\label{eq-pkNsum}
		p_k(N)=\sum_{a=0}^N p(a)\,p_{k-1}(N-a).
	\end{equation}
	Equivalently, the generating function is
	\begin{align*}
		\sum_{N=0}^{\infty}p_k(N)z^N
		=\prod_{n=1}^{\infty}(1-z^n)^{-k}.
	\end{align*}
	General methods for obtaining asymptotic formulas from
	generating functions are discussed in
	\cite{andrews,meinardus, murty}.
	In \cite{barman3}, we established the following
	asymptotic formula for the number of
	$k$-multipartitions of $N$ when $k$ is fixed.
	
	\begin{theorem}[{\cite[Corollary~1.4]{barman3}; see also \cite[Theorem~4]{murty}}]
		\label{thm-murty}
		Let $k$ be any fixed positive integer. Then,
		\begin{equation*}
			p_k(N)=
			\left(\frac{k}{24}\right)^{\frac{k+1}{4}}
			\frac{\exp\left(\frac{2\pi}{\sqrt{6}}\sqrt{Nk}\right)}
			{\sqrt{2}N^{\frac{k+3}{4}}}
			\left(1+O\left(N^{-\frac{1}{2}}\right)\right).
		\end{equation*}
		For $k=1$ this is the Hardy--Ramanujan \cite{hardy} formula $p(N)\sim\frac1{4\sqrt3\,N}\exp\left(\frac{2\pi}{\sqrt6}\sqrt {N}\right)$.
	\end{theorem}
	
	Let $\mathcal{P}_k(N)$ denote the set of all $k$-multipartitions
	of $N$, and let $p_k(N)=|\mathcal{P}_k(N)|$.
	We now state a quantitative result from which
	Theorem~\ref{thm-main} follows.
	
	\begin{theorem}\label{thm-main1}
		Let $G$ be a finite group with an integer-valued character
		table and $k$ conjugacy classes. Let $N$ be sufficiently
		large, and let $q=p^r$ be a prime power
		satisfying
		\begin{equation}\label{eq-qrange}
			q\leq 10^{-3}\,
			\frac{\log(N/k)}{(\log\log N)^2}.
		\end{equation}
		Then the number of pairs
		$(\lambda,\mu)\in \mathcal{P}_k(N)\times\mathcal{P}_k(N)$ such that
		$q\nmid\chi^\lambda_\mu$ is
		\begin{align*}
			O_k\left(
			p_k(N)^2\exp\left(-(\log\log N)^2\right)
			\right).
		\end{align*}
	\end{theorem}
	
	Here, $\chi^\lambda_\mu$ denotes the value of the
	irreducible character of $G\wr S_N$ indexed by $\lambda$
	on the conjugacy class indexed by $\mu$.
	When $G$ is the trivial group, we have $k=1$, and
	Theorem~\ref{thm-main1} reduces to the result of
	Peluse and Soundararajan
	\cite[Theorem~1.1]{peluse3}.
	
	By expressing a fixed positive integer as a product of
	prime powers and applying the union bound, we obtain
	the following corollary.
	\begin{corollary}\label{cor-main2}
		Let $G$ be as in Theorem~\ref{thm-main}, and let $K$
		be a fixed positive integer. Then the proportion of
		entries in the character table of $G\wr S_N$ divisible
		by $K$ tends to $1$ as $N\to\infty$.
	\end{corollary}
	
	Our results also have applications to Weyl groups.
	All Weyl groups have rational-valued characters, and
	hence their character values are integers.
	The Weyl groups of types $B_N$ and $C_N$ are both
	isomorphic to
	$\mathbb{Z}/2\mathbb{Z}\wr S_N$.
	Therefore, Theorem~\ref{thm-main} applies directly
	to these two families.
	The remaining infinite family of type $D_N$ can be
	treated using Clifford theory. This gives the following
	corollary.
	
	\begin{corollary}\label{cor-typeDN}
		Let $q$ be a fixed prime power. As $N\to\infty$,
		the proportion of entries in the character table of
		the Weyl group of type $D_N$ that are divisible by
		$q$ tends to $1$. Hence, the same holds
		for every fixed positive integer in place of $q$.
	\end{corollary}

	\section{Acknowledgments}\label{sec-ackn}
	
	J. Barman is grateful to the University Grants Commission
	(UGC), India, for financial support through its Fellowship
	Programme. K. Mahatab is supported by the ARG-MATRICS
	Programme (Grant No.~ANRF/ARGM/2025/002540/MTR).
	
	\section{Preliminaries}
	In this section, we first discuss the conjugacy classes
	and irreducible characters of $G\wr S_N$ and their
	correspondence with multipartitions. We then recall
	some basic definitions concerning hooks, border strips,
	and core partitions. Finally, we state the
	Murnaghan--Nakayama rule for $G\wr S_N$.
	\subsection{Wreath products, conjugacy classes and irreducible characters}
	
	Let $G$ be a finite group. The wreath product of $G$ with
	the symmetric group $S_N$ is defined by
	$G\wr S_N=G^N\rtimes S_N$.
	Its elements are of the form $(g_1,\ldots,g_N;\sigma)$,
	where $g_i\in G$ and $\sigma\in S_N$. The multiplication
	is given by
	\begin{align*}
		&(g_1,\ldots,g_N;\sigma)(h_1,\ldots,h_N;\tau)
		=(g_1h_{\sigma^{-1}(1)},\ldots,
		g_Nh_{\sigma^{-1}(N)};\sigma\tau).
	\end{align*}
	The group $G\wr S_N$ can also be viewed as the group of
	$N\times N$ monomial matrices with nonzero entries in $G$.
	Its order is $|G|^N N!$.
	
	Let $c_1,\ldots,c_k$ denote the conjugacy classes of $G$,
	and let $\chi_1,\ldots,\chi_k$ be its irreducible characters,
	where $\chi_1$ is the trivial character.
	Consider an element $x=(g_1,\ldots,g_N;\sigma)$ of
	$G\wr S_N$. For a cycle $(i_1\,i_2\,\cdots\,i_m)$ of
	$\sigma$, the conjugacy class of
	$g_{i_m}g_{i_{m-1}}\cdots g_{i_1}$ in $G$ is called
	the cycle product associated with this cycle
	\cite[(4.2.1)]{james}.
	
	The conjugacy classes of $G\wr S_N$ can be described
	using $k$-multipartitions of $N$. We recall the following
	result of James and Kerber.
	
	\begin{proposition}[{\cite[Theorem~4.2.8]{james}}]
		The conjugacy classes of $G\wr S_N$ are indexed by
		the elements of $\mathcal{P}_k(N)$. More precisely, the
		conjugacy class of $x$ corresponds to the
		$k$-multipartition
		$\mu=(\mu^{(1)},\ldots,\mu^{(k)})$, where
		$\mu^{(j)}$ is the partition consisting of the
		lengths of those cycles of $\sigma$ whose cycle
		products belong to $c_j$.
	\end{proposition}
	
	Let $V_1,\ldots,V_k$ be the irreducible representations
	of $G$ corresponding to the irreducible characters
	$\chi_1,\ldots,\chi_k$, respectively. We recall the
	following result of James and Kerber.
	\begin{proposition}[{\cite[Theorem~4.4.3]{james}}]
		The irreducible representations of $G\wr S_N$ are
		indexed by the elements of $\mathcal{P}_k(N)$.
		Let $\lambda=(\lambda^{(1)},\ldots,\lambda^{(k)})
		\in\mathcal{P}_k(N)$, and put $a_i=|\lambda^{(i)}|$.
		Then the corresponding irreducible representation
		is given by
		\begin{align*}
			V^\lambda
			=\mathrm{Ind}_{(G\wr S_{a_1})\times\cdots
				\times(G\wr S_{a_k})}^{G\wr S_N}
			\left(
			\boxtimes_{i=1}^{k}
			\left(S^{\lambda^{(i)}}\otimes V_i^{\otimes a_i}\right)
			\right),
		\end{align*}
		where $S^{\lambda^{(i)}}$ denotes the Specht module
		of $S_{a_i}$ corresponding to the partition
		$\lambda^{(i)}$.
	\end{proposition}
	
	Let $\chi^\lambda$ denote the character of
	$V^\lambda$, and let $\chi^\lambda_\mu$ be its
	value on the conjugacy class indexed by $\mu$.
	If $G$ has an integer-valued character table,
	that is, $\chi_i(c_j)\in\mathbb{Z}$ for all
	$1\leq i,j\leq k$, then the character values
	of $G\wr S_N$ are also integers. In particular,
	$\chi^\lambda_\mu\in\mathbb{Z}$ for all
	$\lambda,\mu\in\mathcal{P}_k(N)$
	\cite[Corollary~4.4.11]{james}.
	
	Moreover, by \cite[\S13.1]{serre}, the character
	table of $G$ is integer-valued if and only if
	every element $g\in G$ is conjugate to $g^a$
	for every integer $a$ relatively prime to
	the order of $g$.
	
	Throughout the paper, we assume that $G$ is a finite group
	with an integer-valued character table and $k$ conjugacy classes.
	
	\subsection{Hooks, border strips and the Murnaghan--Nakayama rule}
	
	Let $\alpha$ be a partition. The hook corresponding
	to a box $b$ in the Young diagram of $\alpha$ consists of
	$b$ and all boxes to its right and below it.
	The hook length $\ell(h)$ is the total number of
	boxes in the hook, while its height $ \mathrm{ht}(h)$
	is one less than the number of rows occupied by the hook.
	
	The border strip, also called a rim hook,
	associated with $h$ is denoted by $\bs(h)$.
	It consists of connected boxes along the south-east
	boundary of the Young diagram, starting from the
	rightmost box of $h$ and ending at its bottommost box.
	The border strip $\bs(h)$ has $\ell(h)$ boxes and
	the same height as $h$. Removing $\bs(h)$ from
	$\alpha$ gives a partition $\alpha\backslash\bs(h)$
	of $|\alpha|-\ell(h)$. Conversely, every removable
	border strip of $\alpha$ corresponds to a unique
	hook of $\alpha$. Here, a border strip is a connected
	skew Young diagram containing no $2\times 2$ square.
	
	For a positive integer $t$, a partition $\alpha$
	is called a $t$-core if none of its hook lengths
	is divisible by $t$.
	
	Let $\lambda=(\lambda^{(1)},\ldots,\lambda^{(k)})
	\in\mathcal{P}_k(N)$. For an index $i$ and a hook $h$
	of $\lambda^{(i)}$, we denote by
	$\lambda\backslash_i h$ the $k$-multipartition
	of $N-\ell(h)$ obtained by removing the border
	strip $\bs(h)$ from $\lambda^{(i)}$, while leaving
	the other components unchanged.
	
	Similarly, let $\mu=(\mu^{(1)},\ldots,\mu^{(k)})
	\in\mathcal{P}_k(N)$. If $\mu^{(j)}$ contains a part
	equal to $m$, then $\mu\backslash_j m$ denotes
	the $k$-multipartition obtained by deleting
	one such part from $\mu^{(j)}$.

	\begin{proposition}[Murnaghan--Nakayama rule for $G\wr S_N$]
		Let $\lambda,\mu\in\mathcal{P}_k(N)$, and suppose that
		$\mu^{(j)}$ contains a part equal to $m$. Then
		\begin{equation}\label{eq-MNrule}
			\chi^\lambda_\mu
			=\sum_{i=1}^{k}\chi_i(c_j)
			\sum_{\substack{h\text{ is a hook of }\lambda^{(i)}\\
					\ell(h)=m}}
			(-1)^{\mathrm{ht}(h)}
			\chi^{\lambda\backslash_i h}_{\mu\backslash_j m},
		\end{equation}
		where the inner sum runs over all hooks of length
		$m$ in $\lambda^{(i)}$. Here,
		$\chi^{(\varnothing,\ldots,\varnothing)}
		_{(\varnothing,\ldots,\varnothing)}=1$.
	\end{proposition}
	
	\begin{proof}
		
		
		By the Murnaghan--Nakayama rule for wreath products
		\cite[Theorem~4.4.10]{james}
		(see also \cite[Proposition~2.8]{dong}),
		the character value $\chi^\lambda_\mu$ can be
		expressed as a sum over all rim hook decompositions
		of $\lambda$ corresponding to $\mu$.
		In each decomposition, border strips are removed
		successively from the components of $\lambda$,
		with lengths equal to the parts of $\mu$ in a
		fixed order.
		
		If a border strip is removed from $\lambda^{(i)}$
		corresponding to a part of $\mu^{(j)}$, then it
		contributes the factor $\chi_i(c_j)$ together
		with the sign $(-1)^{\mathrm{ht}(h)}$,
		where $h$ is the corresponding hook.
		
		Since the rule holds for any fixed ordering of
		the parts of $\mu$, we choose the part $m$ of
		$\mu^{(j)}$ first. We then group the rim hook
		decompositions according to the component
		$\lambda^{(i)}$ and the hook $h$ of length $m$
		removed in the first step.
		
		After removing this border strip, the remaining
		decomposition corresponds to
		$\lambda\backslash_i h$ and $\mu\backslash_j m$.
		Summing over all possible components $i$ and
		hooks $h$ of length $m$, we obtain
		\eqref{eq-MNrule}.
	\end{proof}
	
	When $G$ is the trivial group, we have $k=1$,
	and \eqref{eq-MNrule} reduces to the classical
	Murnaghan--Nakayama rule
	\cite[Theorem~2.4.7]{james}.
	
	\subsection{Permutation characters}
	
	Let $\lambda\in\mathcal{P}_k(N)$, and write
	$a_i=|\lambda^{(i)}|$. Let
	$S_{\lambda^{(i)}}$ be the Young
	subgroup of $ S_{a_i}$ corresponding to $\lambda^{(i)}$,
	and put $G_{\lambda^{(i)}}=G\wr S_{\lambda^{(i)}}$.
	Following \cite[Definition~2.9]{dong},
	we define the permutation module
	\begin{equation*}
		M^\lambda=
		\mathrm{Ind}^{G\wr S_N}_{
			G_{\lambda^{(1)}}\times\cdots\times
			G_{\lambda^{(k)}}}
		\left(\boxtimes_{i=1}^k V_i^{\otimes a_i}\right).
	\end{equation*}
	We denote the value of its character on the
	conjugacy class indexed by $\mu$ by $M^\lambda_\mu$.
	
	\begin{definition}[{\cite[Definition~2.10]{dong}}]
		Let $\lambda,\mu\in\mathcal{P}_k(N)$.
		We treat each row of the components of $\lambda$ and $\mu$
		as a distinct object, even if some rows have the same length.
		A \emph{row decomposition} of $\lambda$
		by $\mu$ is a map $\rho$ from the rows
		of $\mu$ to the rows of $\lambda$ such that,
		for each row $\mathfrak r$ of $\lambda$,
		the sum of the lengths of the rows in
		$\rho^{-1}(\mathfrak r)$ equals the length
		of $\mathfrak r$.
		We denote the set of all such row decompositions
		by $RD(\lambda,\mu)$.
	\end{definition}
	
	The following proposition gives a formula for
	the permutation character values.
	
	\begin{proposition}[{\cite[Proposition~2.11]{dong}}]
		\label{prop-permuchar}
		For $\lambda,\mu\in\mathcal{P}_k(N)$, we have
		\begin{equation*}
			M^\lambda_\mu
			=\sum_{\rho\in RD(\lambda,\mu)}
			\prod_{i=1}^k\prod_{j=1}^k
			\chi_i(c_j)^{a_{ij}(\rho)},
		\end{equation*}
		where $a_{ij}(\rho)$ is the number of rows
		of $\mu^{(j)}$ assigned to rows of
		$\lambda^{(i)}$ by $\rho$.
	\end{proposition}
	
	Recall that the dominance order on partitions
	of the same size is defined by
	$\alpha\succeq\beta$ if
	\begin{equation*}
		\sum_{r=1}^s\alpha_r
		\ge \sum_{r=1}^s\beta_r
	\end{equation*}
	for every $s\ge1$.
	For $k$-multipartitions $\lambda$ and $\eta$
	with equal corresponding component sizes,
	we write $\lambda\succeq\eta$ if
	$\lambda^{(i)}\succeq\eta^{(i)}$
	for every $1\le i\le k$.
	
	\begin{lemma}[{\cite[Lemma~2.13]{dong}}]
		\label{lemma-unitriangu}
		For $\lambda,\eta\in\mathcal{P}_k(N)$, we have
		\begin{equation*}
			[M^\lambda:V^\eta]
			=
			\prod_{i=1}^k K_{\eta^{(i)},\lambda^{(i)}}
		\end{equation*}
		if $|\eta^{(i)}|=|\lambda^{(i)}|$
		for every $i$, and
		$[M^\lambda:V^\eta]=0$ otherwise.
		Here $K_{\alpha,\beta}$ denotes the
		Kostka number. In particular, the matrix $\left([M^\lambda:V^\eta]\right)_{
			\lambda,\eta\in\mathcal{P}_k(N)}$
		is unitriangular with respect to a total
		ordering refining the componentwise
		dominance order, and hence is invertible
		over $\mathbb Z$.
	\end{lemma}
	
	\begin{corollary}\label{cor-Lcombina}
		Every irreducible character $\chi^\lambda$
		of $G\wr S_N$ is an integral linear
		combination of the permutation characters
		$M^\eta$, where $\eta\in\mathcal{P}_k(N)$.
	\end{corollary}
	
	\begin{proof}
		By Lemma~\ref{lemma-unitriangu}, the character
		of $M^\lambda$ satisfies
		\begin{equation*}
			M^\lambda_\mu
			=\sum_{\eta\in\mathcal{P}_k(N)}
			[M^\lambda:V^\eta]\chi^\eta_\mu.
		\end{equation*}
		The matrix of multiplicities is unitriangular
		with integer entries, so its inverse also
		has integer entries. Therefore, every
		irreducible character is an integral linear
		combination of the permutation characters.
	\end{proof}

	\subsection{Counting partitions and multipartitions}
	
	We recall the following results on the number
	of partitions and multipartitions.
	
	\begin{lemma}\label{lema-pkNuperb}
		Let $k\ge1$ be fixed. Then there exist
		positive constants $C_k,c_k$ and $N_0(k)$
		such that, for every $N\ge N_0(k)$,
		the following hold.
		\begin{enumerate}
			\item If $0\le t\le \frac{N}{2}$, then
			\begin{align*}
				p_k(N-t)
				\le C_k p_k(N)
				\exp\left(-\frac{\pi\sqrt{k}\,t}{\sqrt{6N}}\right).
			\end{align*}
			
			\item If $\frac{N}{2}<t\le N$, then
			\begin{align*}
				p_k(N-t)
				\le C_k p_k(N)\exp\left(-c_k\sqrt N\right).
			\end{align*}
		\end{enumerate}
	\end{lemma}
	
	\begin{proof}
		$(1)$ Let $N$ be sufficiently large and
		$0\le t\le\frac{N}{2}$. Since $N-t\ge \frac{N}{2}$,
		Lemma~\ref{thm-murty} gives
		\begin{align*}
			\frac{p_k(N-t)}{p_k(N)}
			&\le 2\left(\frac{N}{N-t}\right)^{\frac{k+3}{4}}
			\exp\left(
			-\frac{2\pi\sqrt{k}}{\sqrt6}
			(\sqrt N-\sqrt{N-t})
			\right).
		\end{align*}
		Now, $\frac{N}{N-t}\le2$ and $\sqrt N-\sqrt{N-t}
		=\frac{t}{\sqrt N+\sqrt{N-t}}
		\ge\frac{t}{2\sqrt N}.$
		Therefore,
		\begin{align*}
			p_k(N-t)
			\le 2^{\frac{k+7}{4}}p_k(N)
			\exp\left(
			-\frac{\pi\sqrt{k}\,t}{\sqrt{6N}}
			\right),
		\end{align*}
		which proves $(1)$.
		
		$(2)$ Suppose that $\frac{N}{2}<t\le N$.
		Since $p_k(n)$ is non-decreasing, we have
		\begin{align*}
			p_k(N-t)
			\le p_k\left(N-\left\lfloor\frac N2\right\rfloor\right).
		\end{align*}
		Applying part $(1)$ with
		$t=\lfloor N/2\rfloor$, we obtain
		\begin{align*}
			p_k(N-t)
			&\le C_k p_k(N)
			\exp\left(
			-\frac{\pi\sqrt{k}}{\sqrt{6N}}
			\left\lfloor\frac N2\right\rfloor
			\right)\le C_k p_k(N)\exp(-c_k\sqrt N),
		\end{align*}
		for some constant $c_k>0$ and all
		sufficiently large $N$.
		This completes the proof.
	\end{proof}
	
	We next recall a lemma that will be useful in the proof of Lemma~\ref{lemma-noncorecount}.
	\begin{lemma}[{\cite[Lemma~5]{morotti}}]
		\label{lema-morotti}
		For integers $1\le t\le n$, the number
		of partitions of $n$ that are not
		$t$-cores is at most $(t+1)p(n-t).$
	\end{lemma}
	\section{Proofs of the Main Results}
	In this section, we first define the relation $\sim_q$
	on multipartitions.
	\begin{definition}
		Let $q=p^r$, where $p$ is a prime and $r\ge1$,
		and let $\mu,\nu\in\mathcal{P}_k(N)$.
		We say that $\nu$ is obtained from $\mu$
		by a $q$-reduction if there exist
		$j\in\{1,\ldots,k\}$ and $m\ge1$ such that
		$\nu^{(i)}=\mu^{(i)}$ for every $i\ne j$,
		and $\nu^{(j)}$ is obtained from $\mu^{(j)}$
		by replacing $p^r$ parts of size $m$
		with $p^{r-1}$ parts of size $pm$.
		Let $\sim_q$ denote the equivalence relation
		on $\mathcal{P}_k(N)$ generated by these reductions.
	\end{definition}
	We now state a congruence relation for the character values
	of $G\wr S_N$.
	\begin{proposition}
		\label{prop-congruence}
		If $\mu\sim_q\nu$, then
		\begin{align*}
			\chi^\lambda_\mu
			\equiv\chi^\lambda_\nu\pmod{p^r}
		\end{align*}
		for every $\lambda\in\mathcal{P}_k(N)$.
	\end{proposition}
	
	Starting with $\mu\in\mathcal{P}_k(N)$, we repeatedly
	apply $q$-reductions whenever a part occurs
	at least $p^r$ times in some component.
	Each reduction decreases the total number
	of parts, so the process terminates.
	
	The final multipartition does not depend on the
	order in which the reductions are performed.
	Indeed, reductions involving different components
	or different part sizes can be performed in either
	order and give the same result. If two reductions
	involve the same component and the same part size,
	they produce the same multipartition. Since each
	reduction decreases the total number of parts,
	the process eventually stops, and the final
	multipartition is uniquely determined.
	
	We denote the resulting multipartition by
	$\widetilde{\mu}$. By construction,
	each component $\widetilde{\mu}^{(j)}$
	is obtained from $\mu^{(j)}$ by the
	reduction procedure described in
	\cite[Proposition~2.3]{peluse3}.
	Hence, Proposition~\ref{prop-congruence}
	gives
	\begin{equation}\label{eq-mutilde}
		\chi^\lambda_\mu
		\equiv\chi^\lambda_{\widetilde{\mu}}
		\pmod{p^r}
	\end{equation}
	for every $\lambda\in\mathcal{P}_k(N)$.

	
	\begin{theorem}[Divisibility criterion]
		\label{thm-divicriterion}
		Let $m_1,\ldots,m_r$ be distinct positive integers,
		and let $j_1,\ldots,j_r\in\{1,\ldots,k\}$,
		not necessarily distinct.
		Suppose that $\mu\in\mathcal{P}_k(N)$ is such that,
		for each $1\le i\le r$, the component
		$\mu^{(j_i)}$ contains at least $p^{r-1}$
		parts equal to $m_i$. Define
		\begin{align*}
			\mathcal{T}
			=\mathcal{T}(m_1,\ldots,m_r)
			:=\Big\{
			\sum_{i=1}^r k_i m_i:
			0\le k_i\le p^{r-1}
			\text{ for all }i,\ 
			k_i=p^{r-1}\text{ for some }i
			\Big\}.
		\end{align*}
		If every component $\lambda^{(a)}$,
		$1\le a\le k$, is a $t$-core for every
		$t\in\mathcal{T}$, then
		\begin{align*}
			p^r\mid\chi^\lambda_\mu.
		\end{align*}
	\end{theorem}
	When $r=1$, the result gives
	$\chi^\lambda_\mu=0$ if every component
	of $\lambda$ is an $m_1$-core,
	as used in \cite{dong}.

	
	\begin{proposition}\label{prop-muside}
		Let $k\ge1$ be fixed. Then there exist constants
		$N_1(k)$ and $C>0$ such that, for every
		$N\ge N_1(k)$ and every prime power $q=p^r$
		satisfying \eqref{eq-qrange}, the following holds.
		For all but at most
		\begin{align*}
			Ck\,p_k(N)\exp\left(-\left(\frac{N}{k}\right)^{\frac{1}{20q}}\right)
		\end{align*}
		multipartitions $\mu\in\mathcal{P}_k(N)$,
		there exists a $j$ such that
		$|\mu^{(j)}|\ge \frac{N}{k}$ and
		$\widetilde{\mu}^{(j)}$ contains $r$
		distinct parts $m_1,\ldots,m_r$,
		each occurring at least $p^{r-1}$ times,
		with
		\begin{equation}\label{eq-miupper}
			p^{r-1}m_i>
			\left(1+\frac{1}{6q}\right)
			\frac{\sqrt6}{2\pi}
			\sqrt{\frac Nk}\log\frac Nk
			\qquad(1\le i\le r).
		\end{equation}
	\end{proposition}
	
	
	\begin{proposition}\label{prop-lambdaupper}
		Let $k\ge1$ be fixed. Then there exist
		constants $C'_k>0$ and $N_2(k)$ such that,
		for every $N\ge N_2(k)$, every
		$0<\eta\le1$, and every integer
		\begin{align*}
			t\ge T_0(\eta):=
			(1+\eta)\frac{\sqrt6}{2\pi}
			\sqrt{\frac Nk}\log\frac Nk,
		\end{align*}
		the number of multipartitions
		$\lambda\in\mathcal{P}_k(N)$ for which at least
		one component $\lambda^{(a)}$ is not
		a $t$-core is at most
		\begin{align*}
			C'_k\,p_k(N)\frac{\log N}{N^{\frac{\eta}{2}}}.
		\end{align*}
	\end{proposition}
	The preceding results provide the necessary
	tools to prove Theorem~\ref{thm-main1}.
	
	\begin{proof}[\textbf{\boldmath Proof of Theorem~\ref{thm-main1}}]
		Assume Propositions~\ref{prop-congruence},
		\ref{prop-muside}, \ref{prop-lambdaupper},
		and Theorem~\ref{thm-divicriterion}.
		Let $N$ be sufficiently large, and let
		$q=p^r$ satisfy \eqref{eq-qrange}.
		Set $\eta=\frac{1}{6q}$. We divide the pairs
		$(\lambda,\mu)\in \mathcal{P}_k(N)\times\mathcal{P}_k(N)$ for which
		$q\nmid\chi^\lambda_\mu$ into two cases.
		\begin{case}
		\end{case}
		Suppose that $\mu$ belongs to the exceptional
		set in Proposition~\ref{prop-muside}.
		The number of such pairs is at most
		\begin{align*}
			Ck\,p_k(N)^2
			\exp\left(-\left(\frac{N}{k}\right)^{\frac{1}{20q}}\right).
		\end{align*}
		\begin{case}
		\end{case}
		Suppose that $\mu$ does not belong to the
		exceptional set. By Proposition~\ref{prop-muside},
		there exist an index $j$ and distinct positive
		integers $m_1,\ldots,m_r$ satisfying
		\eqref{eq-miupper}.
		Let $\mathcal{T}=\mathcal{T}(m_1,\ldots,m_r)$.
		For every $t\in\mathcal{T}$, we have
		$t\ge p^{r-1}m_i$ for some $i$. Hence,
		by \eqref{eq-miupper},
		\begin{align*}
			t>
			\left(1+\frac{1}{6q}\right)
			\frac{\sqrt6}{2\pi}
			\sqrt{\frac Nk}\log\frac Nk
			=T_0(\eta).
		\end{align*}
		Moreover,
		\begin{align*}
			|\mathcal{T}|\le r(p^{r-1}+1)^{r-1},
		\end{align*}
		since there are $r$ choices for an index $i$
		with $k_i=p^{r-1}$, and each of the remaining
		$r-1$ entries can take any integer value
		from $0$ to $p^{r-1}$. Thus, there are at most
		$r(p^{r-1}+1)^{r-1}$ such $r$-tuples.

		By Proposition~\ref{prop-lambdaupper},
		for each $t\in\mathcal{T}$, the number of
		multipartitions $\lambda$ having at least
		one component that is not a $t$-core is at most
		\begin{align*}
			C'_k p_k(N)\frac{\log N}{N^{\frac{1}{12q}}}.
		\end{align*}
		Therefore, the number of multipartitions $\lambda$
		having at least one component that is not a $t$-core
		for some $t\in\mathcal{T}$ is at most
		\begin{align*}
			C'_k|\mathcal{T}|p_k(N)
			\frac{\log N}{N^{\frac{1}{12q}}}.
		\end{align*}
		For every remaining $\lambda$, all components
		are $t$-cores for every $t\in\mathcal{T}$.
		Since $\widetilde{\mu}^{(j)}$ contains at least
		$p^{r-1}$ parts equal to each $m_i$,
		Theorem~\ref{thm-divicriterion} gives
		\begin{align*}
			p^r\mid\chi^\lambda_{\widetilde{\mu}}.
		\end{align*}
		By \eqref{eq-mutilde}, we also have
		\begin{align*}
			p^r\mid\chi^\lambda_\mu.
		\end{align*}
		Thus, for each non-exceptional $\mu$,
		there are at most
		\begin{align*}
			C'_k\, r(p^{r-1}+1)^{r-1}
			p_k(N)\frac{\log N}{N^{\frac{1}{12q}}}
		\end{align*}
		choices of $\lambda$ for which
		$q\nmid\chi^\lambda_\mu$.
		
		Combining the two cases, the total number
		of pairs $(\lambda,\mu)$ with
		$q\nmid\chi^\lambda_\mu$ is at most
		\begin{align}\label{eq-final1}
			p_k(N)^2\left(
			Ck\exp\left(-\left(\frac{N}{k}\right)^{\frac{1}{20q}}\right)
			+C'_k \,r(p^{r-1}+1)^{r-1}
			\frac{\log N}{N^{\frac{1}{12q}}}
			\right).
		\end{align}
		It remains to estimate these two terms.
		Write $\mathcal{L}=\log\log N$.
		By \eqref{eq-qrange}, we have
		\begin{align*}
			q\le\frac{10^{-3}\log(\frac{N}{k})}{\mathcal{L}^2}.
		\end{align*}
		Hence,
		\begin{align*}
			\left(\frac{N}{k}\right)^{\frac{1}{20q}}
			=\exp\left(\frac{\log \frac{N}{k}}{20q}\right)
			\ge\exp(50\mathcal{L}^2)
			\quad\text{and}\quad
			N^{-\frac{1}{12q}}
			=\exp\left(-\frac{\log N}{12q}\right)
			\le\exp(-83\mathcal{L}^2).
		\end{align*}
		Since $p^{r-1}+1\le p^r=q$ and
		$r\le\frac{\log q}{\log2}\le q$, we obtain
		\begin{align*}
			r(p^{r-1}+1)^{r-1}
			\le q^r\le\exp\left(\frac{(\log q)^2}{\log2}\right).
		\end{align*}
		Also, $q\le\log N$ for sufficiently
		large $N$, so $\log q\le\mathcal{L}$.
		Therefore,
		\begin{align*}
			r(p^{r-1}+1)^{r-1}
			\le\exp(1.5\mathcal{L}^2).
		\end{align*}
		It follows that the second term in
		\eqref{eq-final1} is bounded by
		\begin{align*}
			C'_k\exp\left(
			1.5\mathcal{L}^2+\mathcal{L}
			-83\mathcal{L}^2
			\right)
			\ll_k\exp(-\mathcal{L}^2).
		\end{align*}
		Similarly, the first term is at most
		\begin{align*}
			Ck\exp\left(-\exp(50\mathcal{L}^2)\right)
			\ll_k\exp(-\mathcal{L}^2).
		\end{align*}
		Hence, the number of pairs
		$(\lambda,\mu)\in\mathcal{P}_k(N)\times\mathcal{P}_k(N)$ satisfying
		$q\nmid\chi^\lambda_\mu$ is
		\begin{align*}
			O_k\left(
			p_k(N)^2\exp(-(\log\log N)^2)
			\right).
		\end{align*}
		This completes the proof of
		Theorem~\ref{thm-main1}.
	\end{proof}
	We now turn to the proof of Theorem~\ref{thm-main}.
	\begin{proof}[\textbf{\boldmath Proof of Theorem~\ref{thm-main} }]
		
		For any fixed prime power $q$, the condition
		\eqref{eq-qrange} holds for all sufficiently
		large $N$. By Theorem~\ref{thm-main1}, the proportion of entries
		in the character table that are not divisible by $q$ is
		\begin{align*}
			O_k\left(\exp(-(\log\log N)^2)\right)
			\longrightarrow 0, 
			\qquad \text{as }N\to\infty.
		\end{align*}
		This proves Theorem~\ref{thm-main}.
	\end{proof}
	The preceding result allows us to prove
	the following Corollary.
	\begin{proof}[\textbf{\boldmath Proof of Corollary~\ref{cor-main2}}]
		Now let $K$ be a fixed positive integer
		with prime power factorization
		$K=q_1\cdots q_s$, where the $q_i$ are
		pairwise coprime. An integer is divisible
		by $K$ if and only if it is divisible
		by every $q_i$. Therefore, the proportion
		of character table entries not divisible
		by $K$ is at most
		\begin{align*}
			\sum_{i=1}^s
			O_k\left(\exp(-(\log\log N)^2)\right)
			\longrightarrow 0
			\qquad \text{as }N\to\infty.
		\end{align*}
		This proves Corollary~\ref{cor-main2}.
	\end{proof}
	\section{Proof of Proposition~\ref{prop-congruence}}
	\label{sec-congruence}
	To prove Proposition~\ref{prop-congruence}, we first establish
	some necessary results, which are discussed below.
	We then use these results to prove the proposition.
	\begin{lemma}\label{lemma-Fxpr}
		Let $F\in\mathbb Z[x_1,\ldots,x_n]$, and write
		$F(x^p)=F(x_1^p,\ldots,x_n^p)$.
		Then, for every integer $r\ge1$, we have
		\begin{align*}
			F(x)^{p^r}
			\equiv F(x^p)^{p^{r-1}}
			\pmod{p^r}.
		\end{align*}
	\end{lemma}
	
	\begin{proof}
		First, suppose that $r=1$. Write
		$F(x)=\sum_\alpha a_\alpha x^\alpha$,
		where $a_\alpha\in\mathbb Z$.
		By the multinomial theorem, we have
		\begin{align*}
			F(x)^p
			=\sum_\alpha a_\alpha^p x^{p\alpha}
			+pR(x),
		\end{align*}
		for some $R\in\mathbb Z[x_1,\ldots,x_n]$,
		since every multinomial coefficient
		involving at least two distinct terms
		is divisible by $p$.
		By Fermat's little theorem,
		$a_\alpha^p\equiv a_\alpha\pmod p$.
		Hence,
		\begin{equation}\label{eq-Fxp}
			F(x)^p=F(x^p)+pR_1(x),
		\end{equation}
		for some $R_1\in\mathbb Z[x_1,\ldots,x_n]$.
		
		Now let $r\ge2$. Taking the $p^{r-1}$-th power on both sides of
		\eqref{eq-Fxp}, we obtain
		\begin{align*}
			F(x)^{p^r}
			=F(x^p)^{p^{r-1}}+\sum_{\ell=1}^{p^{r-1}}
			\binom{p^{r-1}}{\ell}p^\ell
			R_1(x)^\ell
			F(x^p)^{p^{r-1}-\ell}.
		\end{align*}
		For $1\le\ell\le p^{r-1}$, we have
		\begin{align*}
			\binom{p^{r-1}}{\ell}
			=\frac{p^{r-1}}{\ell}
			\binom{p^{r-1}-1}{\ell-1}.
		\end{align*}
		Therefore,
		\begin{align*}
			v_p\left(\binom{p^{r-1}}{\ell}\right)
			\ge r-1-v_p(\ell),
		\end{align*}
		where $v_p(\ell)$ denotes the exponent
		of $p$ dividing $\ell$.
		Since $v_p(\ell)\le\ell-1$, it follows that
		\begin{align*}
			v_p\left(
			\binom{p^{r-1}}{\ell}p^\ell
			\right)
			&\ge r-1-v_p(\ell)+\ell \ge r.
		\end{align*}
		Thus, every term in the sum is divisible
		by $p^r$, which proves the lemma.
	\end{proof}
	
	\begin{remark}
		Taking $F(x)=a$, where $a\in\mathbb Z$,
		Lemma~\ref{lemma-Fxpr} gives
		\begin{align*}
			a^{p^r}\equiv a^{p^{r-1}}\pmod{p^r}.
		\end{align*}
		When $p\nmid a$, this follows from Euler's
		theorem, since $\varphi(p^r)=p^{r-1}(p-1)$.
		When $p\mid a$, both terms are divisible
		by $p^r$, since $p^{r-1}\ge r$.
		The polynomial congruence in
		Lemma~\ref{lemma-Fxpr} allows us to extend
		this result to the sums arising in the
		row decomposition formula.
	\end{remark}
	In the following subsection, we discuss how
	$M^\lambda_\mu$ can be expressed as a coefficient
	of a polynomial. Related results can be found in
	\cite[Chapter~I, Appendix~B]{macdonald} and \cite{kerber}.
	\subsection{Permutation characters as coefficients}
	Fix $\lambda\in\mathcal{P}_k(N)$. For each row of $\lambda$,
	introduce a variable $x_{i,a}$, where
	$1\le i\le k$ and
	$1\le a\le\ell(\lambda^{(i)})$.
	Define
	\begin{align*}
		x^\lambda
		:=\prod_{i=1}^k
		\prod_{a=1}^{\ell(\lambda^{(i)})}
		x_{i,a}^{\lambda^{(i)}_a}.
	\end{align*}
	For $m\ge1$ and $1\le j\le k$, let
	\begin{align*}
		P_{m,j}(x)
		:=\sum_{i=1}^k\chi_i(c_j)
		\sum_{a=1}^{\ell(\lambda^{(i)})}
		x_{i,a}^m
		\in\mathbb Z[x],
	\end{align*}
	where we use the assumption that
	$\chi_i(c_j)\in\mathbb Z$.
	For any polynomial $H$, let $[x^\lambda]H$
	denote the coefficient of $x^\lambda$ in $H$.
	
	\begin{lemma}\label{lemma-permucoef}
		For every $\mu\in\mathcal{P}_k(N)$, we have
		\begin{align*}
			M^\lambda_\mu
			=[x^\lambda]\prod_{j=1}^k
			\prod_{m\in\mu^{(j)}}P_{m,j}(x),
		\end{align*}
		where the inner product runs over all parts
		of $\mu^{(j)}$, counted with multiplicity.
	\end{lemma}

	\begin{proof}
		We label the rows of $\lambda$ by $(i,a)$, where
		$1\le i\le k$ and
		$1\le a\le \ell(\lambda^{(i)})$.
		Similarly, we label the rows of $\mu$ by $(j,b)$,
		where $1\le j\le k$ and
		$1\le b\le \ell(\mu^{(j)})$.
		By the definition of $P_{m,j}(x)$, we have
		\begin{align*}
			P_{\mu^{(j)}_b,j}(x)
			=\sum_{i=1}^k
			\chi_i(c_j)
			\sum_{a=1}^{\ell(\lambda^{(i)})}
			x_{i,a}^{\mu^{(j)}_b}.
		\end{align*}
		When we expand the product
		\begin{align*}
			\prod_{j=1}^k
			\prod_{m\in\mu^{(j)}}P_{m,j}(x),
		\end{align*}
		we choose one term from each factor.
		This is equivalent to assigning each row of $\mu$
		to a row of $\lambda$.
		Thus, each choice determines a map $\rho$
		from the rows of $\mu$ to the rows of $\lambda$.
		
		For each row $(i,a)$ of $\lambda$, let
		$s_{i,a}(\rho)$ be the sum of the lengths
		of the rows of $\mu$ assigned to $(i,a)$.
		Then the corresponding monomial is
		\begin{align*}
			\prod_{i=1}^k
			\prod_{a=1}^{\ell(\lambda^{(i)})}
			x_{i,a}^{s_{i,a}(\rho)}.
		\end{align*}
		This monomial equals $x^\lambda$ if and only if
		\begin{align*}
			s_{i,a}(\rho)=\lambda^{(i)}_a
		\end{align*}
		for every row $(i,a)$ of $\lambda$.
		By the definition of a row decomposition,
		this holds precisely when
		$\rho\in RD(\lambda,\mu)$.
		For such a row decomposition $\rho$,
		the corresponding coefficient is
		\begin{align*}
			\prod_{i=1}^k\prod_{j=1}^k
			\chi_i(c_j)^{a_{ij}(\rho)},
		\end{align*}
		where $a_{ij}(\rho)$ is the number of rows
		of $\mu^{(j)}$ assigned to rows of
		$\lambda^{(i)}$.
		
		Therefore, the coefficient of $x^\lambda$ is
		\begin{align*}
			[x^\lambda]
			\prod_{j=1}^k
			\prod_{m\in\mu^{(j)}}P_{m,j}(x)
			&=\sum_{\rho\in RD(\lambda,\mu)}
			\prod_{i=1}^k\prod_{j=1}^k
			\chi_i(c_j)^{a_{ij}(\rho)}\\
			&=M^\lambda_\mu,
		\end{align*}
		where the last equality follows from
		Proposition~\ref{prop-permuchar}.
	\end{proof}
	
	\begin{example}
		\setcounter{case}{0}
		Let $G=\mathbb{Z}/2\mathbb{Z}$, with conjugacy classes
		$c_1=\{1\}$ and $c_2=\{-1\}$. The irreducible
		characters of $G$ satisfy
		\begin{align*}
			\chi_1(c_1)=\chi_1(c_2)=\chi_2(c_1)=1,
			\qquad \chi_2(c_2)=-1.
		\end{align*}
		Consider the $2$-multipartitions
		\begin{align*}
			\lambda=((4,2),(3)),\qquad
			\mu=((2,2,1),(3,1))\in\cP_2(9).
		\end{align*}
		We compute $M^\lambda_\mu$ using
		Proposition~\ref{prop-permuchar} and
		Lemma~\ref{lemma-permucoef}.
		
		\noindent\textbf{Using Proposition~\ref{prop-permuchar}.}
		We first compute $M^\lambda_\mu$ using the row decomposition formula.
		The components $\lambda^{(1)}=(4,2)$ and
		$\lambda^{(2)}=(3)$ have rows of lengths
		$4,2$ and $3$, respectively.
		Similarly, $\mu^{(1)}=(2,2,1)$ and
		$\mu^{(2)}=(3,1)$ have rows of lengths
		$2,2,1$ and $3,1$, respectively.
		We treat the two rows of length $2$ as
		distinct rows.
		
		A row decomposition assigns each row of $\mu$
		to a row of $\lambda$ so that the sum of
		the assigned row lengths equals the length
		of the corresponding row of $\lambda$.
		Each row contributes a factor of $1$, except when a row of
		$\mu^{(2)}$ is sent to $\lambda^{(2)}$, in which case
		it contributes $-1$.
		
		We consider two cases depending on where the row of length $3$
		in $\mu^{(2)}$ is assigned.
		\begin{case}  
		\end{case}
		Suppose this row is assigned to the row of
		length $3$ in $\lambda^{(2)}$.
		The remaining rows have lengths $2,2,1,1$
		and must fill the rows of lengths $4,2$
		in $\lambda^{(1)}$. Use the rows $2_1,2_2,1_1,1_2$ (the two rows of length $2$ are distinct).
		This can be done in three ways:
		\begin{align*}
			4=2+1+1,\quad 2=2
			&\qquad \text{(two ways)},\\
			4=2+2,\quad 2=1+1
			&\qquad \text{(one way)}.
		\end{align*}
		\begin{center}
			\begin{tikzpicture}[
				scale=0.85,
				box/.style={
					draw=black,
					thick,
					minimum width=0.8cm,
					minimum height=0.8cm,
					inner sep=0pt
				}
				]
				
				\node at (2,1.0) {\textbf{Possibility 1}};
				
				\foreach \x/\t/\c in {
					0/$2_1$/blue!20,
					1/$2_1$/blue!20,
					2/$1_1$/green!20,
					3/$1_2$/orange!20
				}{
					\node[box,fill=\c] at (\x+0.4,0) {\t};
				}
				
				\foreach \x in {0,1}{
					\node[box,fill=purple!20]
					at (\x+0.4,-0.8) {$2_2$};
				}
				
				\begin{scope}[xshift=5.5cm]
					
					\node at (2,1.0) {\textbf{Possibility 2}};
					
					\foreach \x/\t/\c in {
						0/$2_2$/purple!20,
						1/$2_2$/purple!20,
						2/$1_1$/green!20,
						3/$1_2$/orange!20
					}{
						\node[box,fill=\c] at (\x+0.4,0) {\t};
					}
					
					\foreach \x in {0,1}{
						\node[box,fill=blue!20]
						at (\x+0.4,-0.8) {$2_1$};
					}
					
				\end{scope}
				
				\begin{scope}[xshift=2.75cm,yshift=-3cm]
					
					\node at (2,1.0) {\textbf{Possibility 3}};
					
					\foreach \x/\t/\c in {
						0/$2_1$/blue!20,
						1/$2_1$/blue!20,
						2/$2_2$/purple!20,
						3/$2_2$/purple!20
					}{
						\node[box,fill=\c] at (\x+0.4,0) {\t};
					}
					
					\node[box,fill=green!20]
					at (0.4,-0.8) {$1_1$};
					
					\node[box,fill=orange!20]
					at (1.4,-0.8) {$1_2$};
					
				\end{scope}
				
			\end{tikzpicture}
		\end{center}
		Each decomposition has weight $-1$.
		Thus, the total contribution is $-3$.
		
		\begin{case}
			
		\end{case}
		Suppose the row of length $3$ in $\mu^{(2)}$
		is assigned to the row of length $4$
		in $\lambda^{(1)}$.
		One of the two rows of length $1$ must
		also be assigned to this row.
		The row of length $3$ in $\lambda^{(2)}$
		is then filled by a row of length $2$
		and the remaining row of length $1$.
		The other row of length $2$ fills the
		row of length $2$ in $\lambda^{(1)}$.
		
		There are two choices for the row of length
		$2$ and two choices for the row of length $1$.
		Hence, there are four decompositions.
		Two have weight $1$, and the other two
		have weight $-1$. Their total contribution
		is therefore $0$.
		
		Thus, there are seven row decompositions,
		and Proposition~\ref{prop-permuchar} gives
		\begin{align*}
			M^\lambda_\mu=-3+0=-3.
		\end{align*}
		\noindent\textbf{Using Lemma~\ref{lemma-permucoef}.}
		We now compute $M^\lambda_\mu$ using the coefficient formula.
		Since $\lambda^{(1)}$ has two rows and
		$\lambda^{(2)}$ has one row, we introduce
		the variables $x_{1,1},x_{1,2},x_{2,1}$.
		By the definition of $x^\lambda$, we have
		\begin{align*}
			x^\lambda=x_{1,1}^4x_{1,2}^2x_{2,1}^3.
		\end{align*}
		Using the character values of $G$ in
		the definition of $P_{m,j}(x)$, we obtain
		\begin{align*}
			P_{m,1}(x)
			&=x_{1,1}^m+x_{1,2}^m+x_{2,1}^m,\\
			P_{m,2}(x)
			&=x_{1,1}^m+x_{1,2}^m-x_{2,1}^m.
		\end{align*}
		Since $\mu^{(1)}=(2,2,1)$ and
		$\mu^{(2)}=(3,1)$, the product in
		Lemma~\ref{lemma-permucoef} becomes
		\begin{align*}
			P_{2,1}(x)^2P_{1,1}(x)
			P_{3,2}(x)P_{1,2}(x).
		\end{align*}
		Therefore,
		\begin{align*}
			M^\lambda_\mu
			=[x_{1,1}^4x_{1,2}^2x_{2,1}^3]\,
			P_{2,1}(x)^2P_{1,1}(x)
			P_{3,2}(x)P_{1,2}(x).
		\end{align*}
		Thus, we need to find the coefficient of
		$x_{1,1}^{4}x_{1,2}^{2}x_{2,1}^{3}$ in the polynomial
		$P_{2,1}(x)^2P_{1,1}(x)P_{3,2}(x)P_{1,2}(x)$.
		For simplicity, write
		$x=x_{1,1}$, $y=x_{1,2}$ and
		$z=x_{2,1}$. Then
		\begin{align*}
			M^\lambda_\mu
			=&[x^4y^2z^3]\,
			(x^2+y^2+z^2)^2(x+y+z)(x^3+y^3-z^3)(x+y-z)\\
			=&-3,
		\end{align*}
		which agrees with the result obtained
		from Proposition~\ref{prop-permuchar}.
	\end{example}
	
	We are now ready to prove Proposition~\ref{prop-congruence}.
	\begin{proof}[\textbf{\boldmath Proof of Proposition~\ref{prop-congruence}}]
		It is enough to consider a single $q$-reduction.
		Suppose that $\mu^{(j)}=\xi\cup(m^{p^r})$ and $
		\nu^{(j)}=\xi\cup((pm)^{p^{r-1}}),$
		for some partition $\xi$, and
		$\nu^{(i)}=\mu^{(i)}$ for $i\neq j$.
		Fix $\lambda\in\mathcal{P}_k(N)$, and define
		\begin{align*}
			A(x)
			:=\left(\prod_{\substack{1\le i\le k\\i\neq j}}
			\prod_{m'\in\mu^{(i)}}P_{m',i}(x)\right)
			\left(\prod_{m'\in\xi}P_{m',j}(x)\right).
		\end{align*}
		Thus, $A(x)\in\mathbb Z[x]$ is the product
		corresponding to all parts of $\mu$ except
		the $p^r$ selected parts of size $m$ in
		$\mu^{(j)}$. It is also the product
		corresponding to all parts of $\nu$ except
		the $p^{r-1}$ selected parts of size $pm$
		in $\nu^{(j)}$.
		Since
		\begin{align*}
			P_{pm,j}(x)=P_{m,j}(x^p),
		\end{align*}
		Lemma~\ref{lemma-permucoef} gives
		\begin{align*}
			M^\lambda_\mu-M^\lambda_\nu
			=[x^\lambda]\left(
			A(x)\left(
			P_{m,j}(x)^{p^r}
			-P_{m,j}(x^p)^{p^{r-1}}
			\right)\right).
		\end{align*}
		By Lemma~\ref{lemma-Fxpr}, we have
		\begin{align*}
			P_{m,j}(x)^{p^r}
			-P_{m,j}(x^p)^{p^{r-1}}
			\in p^r\,\mathbb Z[x].
		\end{align*}
		Since $A(x)\in\mathbb Z[x]$, it follows that
		\begin{align*}
			M^\lambda_\mu\equiv M^\lambda_\nu
			\pmod{p^r}.
		\end{align*}
		Finally, by Corollary~\ref{cor-Lcombina},
		every irreducible character $\chi^\lambda$
		can be written as an integral linear
		combination of permutation characters:
		\begin{align*}
			\chi^\lambda
			=\sum_{\eta\in \mathcal{P}_k(N)} d_{\lambda\eta}M^\eta,
			\qquad d_{\lambda\eta}\in\mathbb Z.
		\end{align*}
		Therefore,
		\begin{align*}
			\chi^\lambda_\mu-\chi^\lambda_\nu
			=\sum_{\eta\in \mathcal{P}_k(N)} d_{\lambda\eta}
			\left(M^\eta_\mu-M^\eta_\nu\right)\equiv0\pmod{p^r}.
		\end{align*}
		This proves the proposition.
	\end{proof}

	\section{Proof of Theorem~\ref{thm-divicriterion}}\label{sec-criterion}
	In this section, we extends \cite[Proposition~5.2]{peluse3}
	to wreath products. We then use this proposition to
	prove Theorem~\ref{thm-divicriterion}.
	\begin{proposition}\label{prop-keyprop}
		Let $r,m\ge1$, put $R=p^{r-1}$, and fix
		$j\in\{1,\ldots,k\}$. Suppose that
		$\mu\in\mathcal{P}_k(N)$ and that $\mu^{(j)}$
		contains at least $R$ parts equal to $m$.
		Let $\nu\in\cP_k(N-Rm)$ be obtained from
		$\mu$ by removing $R$ parts equal to $m$
		from $\mu^{(j)}$.
		
		Suppose that $\lambda\in\mathcal{P}_k(N)$ and
		each component $\lambda^{(i)}$ is an
		$Rm$-core. Then
		\begin{align*}
			\chi^\lambda_\mu
			=p\sum_{\lambda'}e_{\lambda'}
			\chi^{\lambda'}_\nu,
		\end{align*}
		where $e_{\lambda'}\in\mathbb Z$.
		The sum runs over all $k$-multipartitions
		$\lambda'$ of $N-Rm$ such that each
		$\lambda'^{(i)}$ is obtained from
		$\lambda^{(i)}$ by successively removing
		$n_i\ge0$ border strips of length $m$,
		where
		\begin{align*}
			n_1+\cdots+n_k=R.
		\end{align*}
	\end{proposition}

	\begin{proof}
		Apply the Murnaghan--Nakayama rule \eqref{eq-MNrule}
		$R$ times, using the $R$ selected parts of size
		$m$ in $\mu^{(j)}$. Since these parts belong
		to the same component, their cycle products
		lie in the same conjugacy class $c_j$.
		We obtain
		\begin{equation}\label{eq-Rfold}
			\chi^\lambda_\mu
			=\sum_{(i_1,h_1),\ldots,(i_R,h_R)}
			\left(
			\prod_{s=1}^R
			\chi_{i_s}(c_j)
			(-1)^{\mathrm{ht}(h_s)}
			\right)
			\chi^{\lambda_R}_\nu.
		\end{equation}
		Here $\lambda_0=\lambda$ and
		$\lambda_s=\lambda_{s-1}\backslash_{i_s}h_s$.
		The sum runs over all sequences in which
		$h_s$ is a hook of length $m$ in the
		component $\lambda_{s-1}^{(i_s)}$.
		
		Now group the terms in \eqref{eq-Rfold}
		according to the final multipartition
		$\lambda'=\lambda_R$. For each $i$, let
		\begin{align*}
			n_i
			=\frac{|\lambda^{(i)}|-|\lambda'^{(i)}|}{m}.
		\end{align*}
		Thus, $n_i$ is the number of hooks removed
		from component $i$, and
		$n_1+\cdots+n_k=R$.
		
		Since removing a hook from one component
		does not affect the other components,
		a sequence leading to $\lambda'$ is
		determined by choosing an order for the
		components and a sequence of hook removals
		within each component.
		There are
		\begin{align*}
			\binom{R}{n_1,\ldots,n_k}
		\end{align*}
		ways to choose the order of the components.
		For each $i$, let
		\begin{align*}
			W^{(i)}
			=W_m^{(n_i)}(\lambda^{(i)},\lambda'^{(i)})
		\end{align*}
		denote the number of sequences of $n_i$
		hooks of length $m$ leading from
		$\lambda^{(i)}$ to $\lambda'^{(i)}$,
		where $W^{(i)}=1$ if $n_i=0$.
		Hence, the total number of sequences
		leading to $\lambda'$ is
		\begin{align*}
			\binom{R}{n_1,\ldots,n_k}
			\prod_{i=1}^k W^{(i)}.
		\end{align*}
		Since all the selected parts belong to
		$\mu^{(j)}$, the character factor is
		\begin{align*}
			\prod_{s=1}^R\chi_{i_s}(c_j)
			=\prod_{i=1}^k\chi_i(c_j)^{n_i},
		\end{align*}
		which depends only on $\lambda'$.
		Also, by \cite[Lemma~5.4]{peluse3},
		the sign of each sequence is
		\begin{align*}
			\prod_{s=1}^R(-1)^{\mathrm{ht}(h_s)}
			=\prod_{i=1}^k
			\epsilon_m(\lambda^{(i)},\lambda'^{(i)}).
		\end{align*}
		Therefore,
		\begin{align*}
			\chi^\lambda_\mu
			&=\sum_{\lambda'}
			\left(
			\prod_{i=1}^k
			\chi_i(c_j)^{n_i}
			\epsilon_m(\lambda^{(i)},\lambda'^{(i)})
			\right)\times
			\binom{R}{n_1,\ldots,n_k}
			\left(\prod_{i=1}^kW^{(i)}\right)
			\chi^{\lambda'}_\nu.
		\end{align*}
		Since $\chi_i(c_j)\in\mathbb Z$,
		it remains to prove that
		\begin{align*}
			p\mid
			\binom{R}{n_1,\ldots,n_k}
			\prod_{i=1}^kW^{(i)}.
		\end{align*}
		First, suppose that at least two of the
		$n_i$ are positive. Then there is an
		index $i_0$ such that
		$0<n_{i_0}<R$. Since
		\begin{align*}
			\binom{R}{n_1,\ldots,n_k}
			=\binom{R}{n_{i_0}}
			\binom{R-n_{i_0}}
			{n_1,\ldots,\widehat{n_{i_0}},\ldots,n_k},
		\end{align*}
		the multinomial coefficient is divisible
		by $\binom{R}{n_{i_0}}$.
		As $R=p^{r-1}$, we have
		\begin{align*}
			p\mid\binom{R}{n_{i_0}}.
		\end{align*}
		Thus, the required divisibility holds.
		Note that this case cannot occur when
		$r=1$, since then $R=1$.
		
		Now suppose that only one of the $n_i$
		is positive. Then $n_{i_0}=R$ for
		some $i_0$, and all other $n_i$ are zero.
		In this case,
		\begin{align*}
			W^{(i_0)}
			=W_m(\lambda^{(i_0)},\lambda'^{(i_0)}).
		\end{align*}
		Since $\lambda^{(i_0)}$ is an $Rm$-core,
		by \cite[Lemma~5.5]{peluse3} gives
		\begin{align*}
			p\mid W^{(i_0)}.
		\end{align*}
		Hence, the required divisibility holds
		in this case as well.
		
		Therefore, every coefficient of
		$\chi^{\lambda'}_\nu$ in the above sum
		is divisible by $p$. This proves the
		proposition.
	\end{proof}
	\begin{remark}
		The assumption that the $R$ cycles of length $m$
		have their cycle products in the same conjugacy
		class $c_j$ ensures that the factor
		$\prod_s\chi_{i_s}(c_j)$ does not depend on
		the order of hook removals from different
		components. This condition holds in our
		application because $q$-reductions do not
		move parts between components. Thus, the
		$p^{r-1}$ equal parts in $\widetilde\mu^{(j)}$
		belong to the same component.
	\end{remark}

	\begin{proof}[\textbf{\boldmath Proof of Theorem~\ref{thm-divicriterion}}]
		Let $R=p^{r-1}$. For $0\le u\le r$, let
		$\mu_u$ be obtained from $\mu$ by removing
		$R$ parts equal to $m_i$ from
		$\mu^{(j_i)}$ for each $i>r-u$.
		Thus,
		\begin{align*}
			\mu_u\in\cP_k\left(
			N-R\sum_{i>r-u}m_i
			\right),
			\qquad \mu_0=\mu.
		\end{align*}
		Since the $m_i$ are distinct, these
		selected parts are distinct.
		For $0\le u\le r$, define
		\begin{align*}
			\mathcal{T}_u
			:=\Big\{
			\sum_{i=1}^{r-u}k_im_i:
			0\le k_i\le R,\;
			k_i=R\text{ for some }i\le r-u
			\Big\}.
		\end{align*}
		In particular, $\mathcal{T}_0=\mathcal{T}$ and
		$\mathcal{T}_r=\varnothing$.
		We prove by induction on $u$ that
		\begin{equation}\label{eq-induc}
			\chi^\lambda_\mu
			=p^u\sum_{\lambda''}
			e^{(u)}_{\lambda''}
			\chi^{\lambda''}_{\mu_u},
			\qquad e^{(u)}_{\lambda''}\in\mathbb Z,
		\end{equation}
		where every component of each multipartition
		$\lambda''$ in the sum is a $t$-core
		for all $t\in\mathcal{T}_u$.
		For $u=0$, this follows directly from
		the assumption on $\lambda$.
		
		Now suppose that \eqref{eq-induc}
		holds for some $u<r$. Let $\lambda''$
		be a multipartition appearing in the sum,
		and put
		\begin{align*}
			m=m_{r-u},\qquad j=j_{r-u}.
		\end{align*}
		Since $Rm\in\mathcal{T}_u$, every component
		of $\lambda''$ is an $Rm$-core.
		Also, $\mu_u^{(j)}$ contains at least
		$R$ parts equal to $m$. Removing
		$R$ of these parts gives $\mu_{u+1}$.
		Therefore, by Proposition~\ref{prop-keyprop},
		\begin{align*}
			\chi^{\lambda''}_{\mu_u}
			=p\sum_{\lambda'''}
			e_{\lambda'''}
			\chi^{\lambda'''}_{\mu_{u+1}},
		\end{align*}
		where each component $\lambda'''^{(a)}$
		is obtained from $\lambda''^{(a)}$
		by successively removing $n_a\le R$
		border strips of length $m$. We now show that every component of
		$\lambda'''$ is a $t$-core for all
		$t\in\mathcal{T}_{u+1}$.
		Fix $t\in\mathcal{T}_{u+1}$. For every
		$0\le k'\le R$, we have
		\begin{align*}
			t+k'm\in\mathcal{T}_u,
		\end{align*}
		because $t$ already has a coefficient
		equal to $R$ at some index
		$i\le r-u-1$. Hence, every component
		$\lambda''^{(a)}$ is a
		$(t+k'm)$-core for $0\le k'\le R$.
		
		By \cite[Lemma~5.3]{peluse3}, removing
		one border strip of length $m$
		preserves the $(t+k'm)$-core property
		for $0\le k'\le R-1$.
		Repeating this argument $n_a$ times
		shows that $\lambda'''^{(a)}$ is a
		$(t+k'm)$-core for
		$0\le k'\le R-n_a$.
		In particular, $\lambda'''^{(a)}$
		is a $t$-core.
		Thus, substituting the above expression
		for $\chi^{\lambda''}_{\mu_u}$ in
		\eqref{eq-induc}, we obtain
		\eqref{eq-induc} with $u+1$
		in place of $u$.
		
		Finally, taking $u=r$ in
		\eqref{eq-induc}, we obtain
		\begin{align*}
			\chi^\lambda_\mu
			=p^r\sum_{\lambda''}
			e^{(r)}_{\lambda''}
			\chi^{\lambda''}_{\mu_r}.
		\end{align*}
		Since the coefficients and character
		values are integers, it follows that $p^r\mid\chi^\lambda_\mu.$
		This completes the proof.
	\end{proof}
	
	\section{Proof of Proposition~\ref{prop-lambdaupper}}\label{sec-counting}
	Before proving Proposition~\ref{prop-lambdaupper},
	we state and prove the following lemma, which
	will be useful in the proof.
	\begin{lemma}\label{lemma-noncorecount}
		Let $1\le t\le N$ be integers. Then the number
		of $\lambda\in\mathcal{P}_k(N)$ for which at least one
		component is not a $t$-core is at most $k(t+1)p_k(N-t).$
	\end{lemma}
	
	\begin{proof}
		Fix $a\in\{1,\ldots,k\}$.
		For each $n=|\lambda^{(a)}|$, we use
		Lemma~\ref{lema-morotti} to count the
		multipartitions for which $\lambda^{(a)}$
		is not a $t$-core. Their number is at most
		\begin{align*}
			&\sum_{n=t}^{N}
			(t+1)p(n-t)p_{k-1}(N-n)=(t+1)\sum_{n'=0}^{N-t}
			p(n')p_{k-1}(N-t-n')=(t+1)p_k(N-t),
		\end{align*}
		where the last equality follows from
		\eqref{eq-pkNsum}.
		Summing over all $a\in\{1,\ldots,k\}$
		gives the required bound.
	\end{proof}
	
	We now use this lemma to prove
	Proposition~\ref{prop-lambdaupper}.
	\begin{proof}[\textbf{\boldmath Proof of Proposition~\ref{prop-lambdaupper}}]
		\setcounter{case}{0}
		Let $t\ge T_0=T_0(\eta)$. If $t>N$,
		then every component of $\lambda$ is a
		$t$-core, since its size is at most $N$.
		Thus, we may assume that $t\le N$.
		Let $\beta=\frac{\pi\sqrt{k}}{\sqrt{6N}}.$
		\begin{case}\textbf{\boldmath $T_0\le t\le \frac{N}{2}$.} 
		\end{case}
		By Lemma~\ref{lemma-noncorecount} and
		Lemma~\ref{lema-pkNuperb}(1), the number
		of such multipartitions is at most
		\begin{align*}
			k(t+1)C_kp_k(N)\exp(-\beta t).
		\end{align*}
		Let $g(t)=(t+1)\exp(-\beta t).$
		Since $g$ is decreasing for
		$t\ge1/\beta$, and
		\begin{align*}
			T_0
			=(1+\eta)\frac{\sqrt6}{2\pi}
			\sqrt{\frac Nk}\log\frac Nk
			\ge\frac1\beta
		\end{align*}
		for all sufficiently large $N$,
		we have $g(t)\le g(T_0)$.
		Also, $\beta T_0
		=\frac{\eta+1}{2}\log\frac Nk.$
		Therefore,
		\begin{align*}
			g(T_0)
			&=(T_0+1)
			\left(\frac Nk\right)^{-\frac{\eta+1}{2}}\le
			2\frac{\sqrt6}{\pi}
			\sqrt{\frac Nk}\,(\log N)
			\left(\frac Nk\right)^{-\frac{\eta+1}{2}}\le 2\sqrt{k}\,\frac{\log N}{N^{\frac{\eta}{2}}},
		\end{align*}
		where we use $0<\eta\le1$ and
		assume that $N$ is sufficiently large.
		Hence, the number of such multipartitions
		is at most
		\begin{align*}
			2k^{\frac{3}{2}}C_kp_k(N)\frac{\log N}{N^{\frac{\eta}{2}}}.
		\end{align*}
		\begin{case}\textbf{\boldmath $\frac{N}{2}<t\le N$.} 
		\end{case}
		By Lemma~\ref{lemma-noncorecount} and
		Lemma~\ref{lema-pkNuperb}(2), the number
		of such multipartitions is at most
		\begin{align*}
			k(N+1)C_kp_k(N)\exp\left(-c_k\sqrt{N}\right)
			\le p_k(N)\frac{\log N}{N^{\eta/2}},
		\end{align*}
		for all sufficiently large $N$, since
		$0<\eta\le1$. This proves the proposition.
	\end{proof}

	\section{Proof of Proposition~\ref{prop-muside}}\label{sec-propmuside}
	
	We first recall the following result of
	Peluse and Soundararajan for ordinary partitions.
	Its proof uses estimates for partitions into
	powers of $p$ and a saddle-point bound for the
	partition generating function
	(see \cite[Sections~8--9]{peluse3}).
	
	\begin{proposition}[{\cite[Proposition~2.3]{peluse3}}]
		\label{prop-peluse}
		There exist absolute constants $N_0$ and $C$
		such that the following holds.
		Let $N\ge N_0$, and suppose that
		\begin{align*}
			q=p^r\le
			10^{-3}\frac{\log N}{(\log\log N)^2}.
		\end{align*}
		For a partition $\alpha\vdash N$, let
		$\widetilde\alpha$ be the partition obtained
		by replacing $p^r$ equal parts $m$ with
		$p^{r-1}$ equal parts $pm$. We repeat this
		process until no part occurs more than
		$p^r-1$ times.
		Then, for all but at most
		\begin{align*}
			Cp(N)\exp\left(-N^{\frac{1}{20q}}\right)
		\end{align*}
		partitions $\alpha\vdash N$, the partition
		$\widetilde\alpha$ contains $r$ distinct
		part sizes $m_1,\ldots,m_r$, each occurring
		at least $p^{r-1}$ times, such that
		\begin{align*}
			p^{r-1}m_i>
			\left(1+\frac{1}{6q}\right)
			\frac{\sqrt6}{2\pi}\sqrt N\log N,
			\qquad 1\le i\le r.
		\end{align*}
	\end{proposition}
	We call a partition $\alpha\vdash N$
	$q$-exceptional if it belongs to the
	exceptional set in Proposition~\ref{prop-peluse}.
	Let $b_q(N)$ denote the number of
	$q$-exceptional partitions of $N$.
	
	We now use Proposition~\ref{prop-peluse}
	to prove Proposition~\ref{prop-muside}.
	
	\begin{proof}[\textbf{\boldmath Proof of Proposition~\ref{prop-muside}}]
		Let $N$ be sufficiently large so that $\frac Nk\ge\max\{N_0,e^{e^2}\},$
		and suppose that $q$ satisfies
		\eqref{eq-qrange}.
		The function $x\longmapsto
		\frac{\log x}{(\log\log x)^2}$
		is increasing for $x\ge e^{e^2}$.
		Hence, for every integer $a$ with
		$\frac{N}{k}\le a\le N$, we have
		\begin{align*}
			q\le10^{-3}
			\frac{\log \frac{N}{k}}{(\log\log N)^2}\le10^{-3}
			\frac{\log \frac{N}{k}}{\left(\log\log \frac{N}{k}\right)^2}\le10^{-3}
			\frac{\log a}{(\log\log a)^2}.
		\end{align*}
		Thus, Proposition~\ref{prop-peluse}
		applies to partitions of $a$, giving
		\begin{align*}
			b_q(a)\le Cp(a)\exp\left(-a^{\frac{1}{20q}}\right).
		\end{align*}
		Now let $\mu\in\mathcal{P}_k(N)$. Since $\sum_{j=1}^k|\mu^{(j)}|=N,$
		there is at least one index $j$ such that
		$|\mu^{(j)}|\ge \frac{N}{k}$.

		Suppose that at least one such component
		$\mu^{(j)}$ is not exceptional in the sense
		of Proposition~\ref{prop-peluse}.
		Fix this $j$, and put $a=|\mu^{(j)}|$.
		By Proposition~\ref{prop-peluse},
		the reduced partition
		$\widetilde\mu^{(j)}$ contains
		$r$ distinct part sizes
		$m_1,\ldots,m_r$, each occurring
		at least $p^{r-1}$ times, such that
		\begin{align*}
			p^{r-1}m_i
			&>
			\left(1+\frac1{6q}\right)
			\frac{\sqrt6}{2\pi}\sqrt a\log a\ge
			\left(1+\frac1{6q}\right)
			\frac{\sqrt6}{2\pi}
			\sqrt{\frac {N}{k}}\log\frac {N}{k},
		\end{align*}
		where the last inequality follows
		from $a\ge \frac {N}{k}$ and the fact that
		$\sqrt{x}\log x$ is increasing for
		$x\ge \frac {N}{k}$.
		Hence, \eqref{eq-miupper} holds. Therefore, the exceptional multipartitions
		must belong to
		\begin{align*}
			\bigcup_{j=1}^k\mathcal E_j,
		\end{align*}
		where $\mathcal E_j$ is the set of
		$\mu\in\mathcal{P}_k(N)$ such that
		$|\mu^{(j)}|\ge \frac {N}{k}$ and
		$\mu^{(j)}$ is exceptional in
		Proposition~\ref{prop-peluse}.
		For each fixed $j$, counting according
		to $a=|\mu^{(j)}|$, we obtain
		\begin{align*}
			|\mathcal E_j|
			=\sum_{\frac {N}{k}\le a\le N}
			b_q(a)p_{k-1}(N-a)
			&\le C\exp\left(
			-\left(\frac Nk\right)^{\frac{1}{20q}}
			\right)
			\sum_{a=0}^N p(a)p_{k-1}(N-a)\\
			&=Cp_k(N)\exp\left(
			-\left(\frac Nk\right)^{\frac{1}{20q}}
			\right),
		\end{align*}
		where the last equality follows
		from \eqref{eq-pkNsum}.
		
		Summing over $j=1,\ldots,k$, we conclude that
		the number of exceptional multipartitions
		is at most
		\begin{align*}
			Ck\,p_k(N)\exp\left(
			-\left(\frac Nk\right)^{\frac{1}{20q}}
			\right).
		\end{align*}
		This proves the proposition.
		
	\end{proof}

	\section{Weyl groups of type \texorpdfstring{$D_N$}{D}}
	\label{sec-typeD}
	
	In this section, we prove Corollary~\ref{cor-typeDN}.
	Recall that the Weyl group of type $D_N$
	is an index-$2$ subgroup of the Weyl group
	of type
	$B_N=(\mathbb{Z}/2\mathbb{Z})\wr S_N$.
	Here, $D_N$ denotes the Weyl group,
	not the dihedral group.
	The conjugacy classes of $B_N$ are indexed
	by $2$-multipartitions
	$\lambda=(\lambda^{(1)},\lambda^{(2)})$
	of $N$, where $|\lambda^{(1)}|+|\lambda^{(2)}|=N.$

	\begin{lemma}[{\cite[Proposition~3.4.12]{geck}}]
		\label{lemma-Geck1}
		A conjugacy class of $B_N$ is contained
		in $D_N$ if and only if the corresponding
		$2$-multipartition
		$(\lambda^{(1)},\lambda^{(2)})$
		has an even number of parts in
		$\lambda^{(2)}$.
		Such a conjugacy class remains a single
		conjugacy class of $D_N$, except when
		$\lambda^{(2)}=\varnothing$ and all
		parts of $\lambda^{(1)}$ are even.
		In this case, it splits into two
		conjugacy classes of $D_N$.
	\end{lemma}
	
	Since $D_N$ is a normal subgroup of
	$B_N$ of index $2$, Clifford theory
	gives the following description of
	its irreducible characters.
	
	\begin{lemma}[{\cite[\S5.6.1]{geck}}]
		\label{lemma-Geck2}
		Let
		$\lambda=(\lambda^{(1)},\lambda^{(2)})$
		and
		$\lambda'=(\lambda^{(2)},\lambda^{(1)})$
		be $2$-multipartitions of $N$.
		Then the following hold.
		\begin{enumerate}
			\item
			If $\lambda^{(1)}\neq\lambda^{(2)}$,
			then $\chi^\lambda\!\downarrow_{D_N}$
			is irreducible, and
			\begin{align*}
				\chi^\lambda\!\downarrow_{D_N}
				=
				\chi^{\lambda'}\!\downarrow_{D_N}.
			\end{align*}
			Distinct unordered pairs
			$\{\lambda,\lambda'\}$ give distinct
			irreducible characters of $D_N$.
			
			\item
			If $\lambda^{(1)}=\lambda^{(2)}$,
			then $N$ is even, and
			\begin{align*}
				\chi^\lambda\!\downarrow_{D_N}
				=
				\chi^{\lambda,+}+\chi^{\lambda,-},
			\end{align*}
			where $\chi^{\lambda,+}$ and
			$\chi^{\lambda,-}$ are distinct
			irreducible characters of $D_N$.
		\end{enumerate}
		
		Every irreducible character of $D_N$
		arises from one of these two cases.
	\end{lemma}
	
	It follows from Lemma~\ref{lemma-Geck2}
	that the number of irreducible characters
	of $D_N$, which is also the number of
	its conjugacy classes, is
	\begin{align*}
		k(D_N)
		:=|\mathrm{Irr}(D_N)|
		=\frac{p_2(N)-d(N)}{2}+2d(N)
		=\frac{p_2(N)+3d(N)}{2},
	\end{align*}
	where
	\begin{align*}
		d(N)=
		\begin{cases}
			p(N/2), & \text{if $N$ is even},\\
			0, & \text{if $N$ is odd}.
		\end{cases}
	\end{align*}
	
	We are now ready to prove
	Corollary~\ref{cor-typeDN}.

	\begin{proof}[\textbf{\boldmath Proof of Corollary~\ref{cor-typeDN}}]
		The number of irreducible characters of
		$D_N$ obtained by restricting
		$\chi^\lambda$ with
		$\lambda^{(1)}\neq\lambda^{(2)}$ is $\frac{p_2(N)-d(N)}{2}.$
		By Theorem~\ref{thm-murty}, we have $\frac{d(N)}{p_2(N)}\longrightarrow0$ as $N\to\infty.$
		Thus, the proportion of such irreducible
		characters of $D_N$ tends to $1$.
		
		Now, by Lemma~\ref{lemma-Geck1}, the
		conjugacy classes of $B_N$ contained in
		$D_N$ remain single conjugacy classes,
		except for $d(N)$ classes, each of which
		splits into two classes of $D_N$. Since $D_N$ is a subgroup of $B_N$ of
		index $2$, at least half of the conjugacy
		classes of $B_N$ are contained in $D_N$.
		Hence, these classes form a positive
		proportion of the conjugacy classes of $B_N$.
		
		By Corollary~\ref{cor-main2}, applied
		with $G=\mathbb{Z}/2\mathbb{Z}$, the
		proportion of entries in the character
		table of $B_N$ not divisible by any
		fixed positive integer $K$ tends to $0$.
		Therefore, the same holds when we
		restrict to the columns corresponding
		to conjugacy classes contained in $D_N$.
		Moreover, the number of split classes
		and the number of irreducible characters
		arising from split restrictions are both
		negligible compared with $p_2(N)$.
		
		It follows that the proportion of entries
		in the character table of $D_N$ not
		divisible by $K$ tends to $0$.
		Hence, the proportion of entries divisible
		by $K$ tends to $1$ as $N\to\infty$.
	\end{proof}

	\subsection*{AI Declaration}
	The authors declare that no artificial intelligence tools or AI-assisted technologies were
	used in the preparation, writing, or mathematical analysis of this manuscript.

\end{document}